# Optimization of the first mixed boundary value problem for parabolic differential inclusions in a three spatial dimension

**Elimhan N. Mahmudov**

*Department of Mathematics, Istanbul Technical University, Istanbul, Turkey,*

*Azerbaijan National Academy of Sciences Institute of Control Systems, Baku, Azerbaijan.*

*elimhan22@yahoo.com, ORCID ID: 0000-0003-2879-6154*

**Abstract.** The paper is devoted to the optimization of a first mixed boundary value problem for parabolic differential inclusions (DFIs) with Laplace operator. For this, a problem with a parabolic discrete inclusion is defined, which is the main auxiliary problem. With the help of locally adjoint mappings, necessary and sufficient conditions for the optimality of parabolic discrete inclusions are proved. Then, using the method of discretization of parabolic DFIs and the already obtained optimality conditions for discrete inclusions, the necessary and sufficient conditions for the discrete-approximate problem are formulated in the form of the Euler-Lagrange type inclusion. Thus, using specially proved equivalence theorems, without which it would hardly be possible to obtain the desired result for the problem posed, we establish sufficient optimality conditions for a parabolic DFIs. To demonstrate the above approach, some linear problems and polyhedral optimization with inclusions of parabolic type are investigated.

1. INTRODUCTION

Over the past decades, great progress has been made in various areas of optimal control problems described by ordinary [1,6,8,10,11,13,14,16,19,20,24,28-31,36-39] and partial differential equations/inclusions [4,5,15,21-23,26]. In the paper [12] the averaging method is used to study singularly perturbed differential inclusions in an evolutionary triplet. The upper and lower limits of the set of slow solutions are estimated when the small parameter tends to zero and there is a system with an upper or lower averaging limit. Sufficient conditions for the existence of the limit of the set of solutions are given. The report [40] is devoted to second order discrete approximations to differential inclusions. The approximations are of the form of discrete inclusions with right-hand sides, which are explicitly described for some classes of differential inclusions. The paper [18] investigates some aspects of the form of solutions of linear parabolic partial differential equations in one spatial dimension. The work [3] investigates boundary value problems for systems of Hamilton-Jacobi-Bellman first-order partial differential equations and variational inequalities, the solutions of which obey the viability constraints. They are motivated by some management problems and financial mathematics. The existence and uniqueness of such solutions in the class of closed multivalued mappings is proved. In the paper [34], problems of optimal control of dynamical systems controlled by functional differential inclusions of neutral type with constraints are considered. Developing the method of discrete approximations and using advanced tools of generalized differentiation, a variational analysis of neutral functional-differential inclusions was carried out and new necessary optimality conditions were obtained, both of the Euler-Lagrange and Hamiltonian types. In the paper [9] by using the Baire category method are proven an existence result for boundary value problem of Dirichlet type, for non-convex ordinary differential inclusions under Caratheodory assumptions. In [17] PDIs on the basis of the diffusion properties of weak solutions of stochastic partial DFIs, some existence theorems and some properties of solutions are given. In the book [2] in the framework of the theory of viability, first order PDIs, the solutions of which are feedbacks, are constructed and investigated. In [7] are studied an optimal control problem given by hyperbolic type DFIs with boundary conditions and endpoint constraints. An approach concerning second-order optimality conditions is proposed. The related paper [27] is devoted to the optimal control problem defined by hyperbolic discrete and differential inclusions of Darboux type. The results are extended to non-convex problems. An approach concerning necessary and sufficient conditions for optimality is proposed. In order to formulate sufficient conditions of optimality for problem with partial DFIs the approximation method is used.

The present paper is dedicated to one of the difficult and interesting fields, where the main goal is to derive sufficient optimality conditions for problems with parabolic DFIs. Here we obtain the results in the same spirit as in [27], where some control problems controlled by Darboux type DFIs are considered. A classic example of a propagation problem defined by a parabolic DFI is nonstationary diffusion of heat in a solid. In general, for parabolic DFIs, diffusion, Brownian motion, and flux of heat or electric charges all provide good interpretations. The paper is organized in the following order:

In *Section* 2, the needed concepts and results, such as, Hamiltonian function $H_F$ and argmaximum sets of a set-valued mapping $F$, the locally adjoint mappings (LAMs), etc., from the book of Mahmudov [25], are given and the problem with discrete and continuous parabolic inclusions are stated. Moreover, some "non-degeneracy condition", that is, the standard condition for convex analysis about the existence of an interior point is formulated.

*In Section 3,* for a problem with discrete inclusions of parabolic type, we prove a theorem on necessary and sufficient conditions for optimality. It is shown that under non-degeneracy condition the necessary conditions are also sufficient for the optimality. In the proof of this theorem, the main role is played by the reduction of the problem posed to the problem of mathematical programming and, as a consequence, the application of its methods.

*In Section 4* we introduce difference operators defined on three-point models, and using difference approximations of the Laplace operator and grid functions on a uniform grid, we approximate the problem with parabolic PFIs. In what follows, with a skilful definition of an auxiliary set-valued mapping $G(\cdot, x, y, t): \mathbb{R}^{5n} \rightrightarrows \mathbb{R}^n$ involving the original set-valued mapping $F(\cdot, x, y, t): \mathbb{R}^n \rightrightarrows \mathbb{R}^n$ we reduce this problem to a problem with a discrete inclusion and apply the results obtained in Section 3 and obtain a necessary and sufficient optimality condition for a discrete approximate problem. But, in turn, the latter is possible due to the LAM equivalence theorem, providing the transition from LAM $G^*$ to LAM $F^*$, without which it would hardly be possible to obtain the desired result for problem with parabolic DFIs (PC). For this, auxiliary results are obtained on the equivalence of the cone of tangent directions and the argmaximum sets.

*In Section* 5, we use the results of Section 4 to obtain sufficient optimality conditions for parabolic PDIs with Laplace operator. The formulation of sufficient conditions is realized by going over to the formal limit in the discrete-approximate problem, as discrete steps tend to zero. At the end of section is considered some linear problem with parabolic DFIs. To this end, it is likely that, by intuition, the Euler-Lagrange adjoint DFIs for the continuous problem (PC) can be established, but this intuition by the same time can be accompanied by many unsuccessful attempts. Contrary to this intuition, we establish such a result for parabolic PDEs. Our methodology with a transition to a discrete, and then a discrete-approximate problem guarantees to construct the parabolic type Euler-Lagrange adjoint DFIs and the "endpoint" and boundary conditions for the problem with parabolic DFIs step by step. Without any doubt, this is one of the advantages of this paper.

*Section* 6 is devoted to the problem with parabolic polyhedral DFIs. Here, the integrand of the Lagrange functional is also a polyhedral function. The Euler-Lagrange type inclusion is constructed, and necessary and sufficient optimality condition for problem with parabolic polyhedral DFIs involving Laplace operator are proved. Apparently, such difference problems, in addition to being of independent interest, can play an important role also in computational procedures.

## 2. NECEESARY FACTS AND PROBLEM STATEMENTS

Although all the necessary definitions and concepts can be found in the book [25], for the convenience of readers, they are given in this section. Let $(u,v)$ and $\langle u,v \rangle$ be a pair and inner product of elements $u,v$ in $n$-dimensional Euclidean space $\mathbb{R}^n$, respectively. A set-valued mapping $F: \mathbb{R}^n \rightrightarrows \mathbb{R}^n$ is convex if its graph $\operatorname{gph} F = \{(u,v): v \in F(u)\}$ is convex in $\mathbb{R}^{2n}$, $F$ is convex valued if $F(u)$ is a convex set for each $u \in \operatorname{dom} F = \{u: F(u) \neq \varnothing\}$. For a set-valued mapping, we introduce the Hamiltonian function and the notation of the argmaximum set:

$$H_F(u,v^*) = \sup_v \{\langle v,v^* \rangle: v \in F(u)\}, v^* \in \mathbb{R}^n; \quad F_A(u;v^*) = \{v \in F(u): \langle v,v^* \rangle = H_F(u,v^*)\}.$$

For a convex set-valued $F$ we let $H_F(u,v^*) = +\infty$, if $F(u) = \varnothing$.

The cone of tangent directions at the point $(u,v) \in \operatorname{gph} F$ will be denoted by $K_F(u,v)$ $\equiv K_{\operatorname{gph} F}(u,v)$ and for a convex set-valued mapping $F$ it is defined at a point $(u,v) \in \operatorname{gph} F$ as follows

$$K_{\operatorname{gph} F}(u,v) = \{\operatorname{cone}(\operatorname{gph} F - (u,v)) = \{(\bar{u}, \bar{v}):$$
$$\bar{u} = \mu(\tilde{u} - u), \bar{v} = \mu(\tilde{v} - v), \mu > 0, (\tilde{u}, \tilde{v}) \in \operatorname{gph} F \}.$$

Recall that a set-valued mapping $F^*(\cdot\,;(u,v)): \mathbb{R}^n \rightrightarrows \mathbb{R}^n$ defined by

$$F^*(v^*;(u,v)) = \{u^* : (u^*, -v^*) \in K_F^*(u,v)\}$$

is called the locally adjoint mapping (LAM) to $F$ at a point $(u,v)$ where $K_F^*(u,v)$ $\equiv K^*_{\operatorname{gph} F}(u,v)$ is the cone dual to the cone $K_F(u,v)$. In what follows, we mainly use the equivalent definition of LAM in terms of the Hamiltonian function A set-valued mapping is called the LAM to "nonconvex" mapping $F$ at a point $(u,v) \in \operatorname{gph} F$:

$$F^*(v^*;(u,v)) := \{u^* : H_F(\tilde{u}, v^*) - H_F(u, v^*) \leq \langle u^*, \tilde{u} - u \rangle, \forall \tilde{u} \in \mathbb{R}^n\}, v \in F_A(u; v^*).$$

Since for a convex set-valued mapping $H_F(\cdot, v^*)$ is concave function, by Theorem 2.1 [25, p.62] the following assertion holds $F^*(v^*;(u,v)) = \partial H_F(u,v^*)$, $v \in F_A(u,v^*)$, where $\partial H_F(u,v^*)$ $= -\partial[-H_F(u,v^*)]$. In fact, the given in the paper notion LAM is closely related to the coderivative concept of Mordukhovich [33-35], which is essentially different for nonconvex mappings. The notion of coderivative has been introduced for set-valued mappings in terms of the basic normal cone to their graphs and is defined as $D^*F(x,v)(v^*) = \{x^* : (x^*, -v^*)$ $\in N((x,v); \operatorname{gph} F)\}$, where $N((x,v); \operatorname{gph} F)$ is a normal cone at $(x,v) \in \operatorname{gph} F$. In the most interesting settings for the theory and applications, coderivatives are nonconvex-valued and hence are not tangentially /derivatively generated. This is the case of the first coderivative for general finite dimensional set-valued mappings for the purpose of applications to optimal control. The main advantage of the definition of LAM is its simplicity.

The cone $K_A(z)$ of tangent directions of the set $A$ at a point $z=(u,v) \in A \subset \mathbb{R}^{2n}$ is called a local tent [25, p.120] if for each $\bar{z} \in ri\, K_A(z)$ there exists a convex cone $K \subseteq K_A(z)$ and a continuous mapping $\psi(\bar{z})$ defined in a neighborhood of the origin such that

(1) $\bar{z} \in ri\, K$, $Lin K = Lin\, K_A(z)$, where $Lin\, K$ is the linear span of $K$,

(2) $\psi(\bar{z}) = \bar{z} + r(\bar{z})$, $\quad \|\bar{z}\|^{-1} r(\bar{z}) \to 0$, as $\bar{z} \to 0$

(3) $z + \psi(\bar{z}) \in A$, $\bar{z} \in K \cap S_\varepsilon(0)$ for some $\varepsilon > 0$, where $S_\varepsilon(0)$ is the ball of radius $\varepsilon$ and with center the origin.

A function $g$ is said to be proper if it does not take the value $-\infty$ and is not identically equal to $+\infty$.

First, consider the following optimization problem with a discrete analogue of the so-called first mixed parabolic DFIs, labelled by (PD):

$$\text{minimize} \sum_{\substack{(x,y,t) \in L_0 \times S_0 \times T_0 \\ (x,y) \neq (0,0), (L,S)}} g(u_{x,y,t}, x, y, t), \tag{1}$$

(PD) $\quad u_{x+1,y,t} \in Q(u_{x-1,y,t}, u_{x,y-1,t}, u_{x,y,t}, u_{x,y+1,t}, u_{x,y,t+1}, x, y, t)$, $(x,y,t) \in L_1 \times S_1 \times T_1$, (2)

$$u_{x,y,0} = \alpha_{xy}, \; u_{x,0,t} = \beta_{x0t}, \; u_{x,S,t} = \beta_{xSt}, \; u_{0,y,t} = \gamma_{0yt}, \; u_{L,y,t} = \gamma_{Lyt}, \tag{3}$$

where $L, S, T$ are fixed natural numbers, $L_i = \{i, ..., L-i\}$; $S_i = \{i, ..., S-i\}$; $T_i = \{0, ..., T-i\}$ $(i=0,1)$ and $g(\cdot, x, y, t): \mathbb{R}^n \to \mathbb{R}^1 \cup \{+\infty\}$, $Q(\cdot, x, y, t): \mathbb{R}^{5n} \rightrightarrows \mathbb{R}^n$ is a set-valued mapping, $\alpha_{xy}, \beta_{x0t}, \beta_{xSt}, \gamma_{0yt}, \gamma_{Lyt}$ are fixed vectors, for all $x, y, t$, respectively. A set $\{u_{x,y,t}\}$ $= \{u_{x,y,t} : (x,y,t) \in L_0 \times S_0 \times T_0, (x,y) \neq (0,0), (L,S)\}$, is called a feasible solution for the problem (1)-(3) if it satisfies the inclusion (2) and initial and boundary conditions (3). The problem consists in finding a solution $\{\tilde{u}_{x,y,t}\}$ of a parabolic discrete inclusions (2) that minimizes (1).

Section 5 of the present paper is devoted to the study of a first mixed boundary value problem for parabolic DFIs in a three spatial dimension, containing Laplace operator in a bounded region:

$$\text{minimize} \quad J[u(\cdot,\cdot,\cdot)] = \int_0^T \iint_D g(u(x,y,t), x, y, t)\, dx dy dt, \tag{4}$$

(PC) $\quad \dfrac{\partial u(x,y,t)}{\partial t} - \Delta u(x,y,t) \in F(u(x,y,t), x, y, t)$, $(x,y,t) \in D \times [0,T]$, (5)

$$u(x,y,0) = \alpha(x,y), \; u(x,0,t) = \beta_0(x,t), \; u(x,S,t) = \beta_S(x,t),$$
$$u(0,y,t) = \gamma_0(y,t), \; u(L,y,t) = \gamma_L(y,t), \; D = [0,L] \times [0,S], \tag{6}$$

where $F(\cdot, x, y, t): \mathbb{R}^n \rightrightarrows \mathbb{R}^n$ is a convex set-valued mapping, $g(\cdot, x, y, t)$ is proper convex function, $\Delta$ is Laplace's operator: $\Delta = \frac{\partial^2}{\partial x^2} + \frac{\partial^2}{\partial y^2}$ and $\alpha$ and $\beta_0, \beta_S, \gamma_0, \gamma_L$ are given continuous functions, $\alpha: D \to \mathbb{R}^n$, $\beta_0, \beta_S : [0,L] \times [0,T] \to \mathbb{R}^n$; $\gamma_0, \gamma_L : [0,S] \times [0,T] \to \mathbb{R}^n$; $L, S, T$ are real positive numbers. Since parabolic inclusion (5) is first order in time, values of $u(x,t)$ must be specified along the initial time boundary. On the other hand, since equation (5) is the second order in space, the values of $u(x,t)$ must be specified along boundaries of the region $D \times [0,T]$. We label this problem (PC). The problem consists in finding a solution $\tilde{u}(x,y,t)$ of a first mixed boundary value problem (PC) that minimizes (4). For convenience, we assume throughout the context that feasible solutions are classical solutions; let $P_T = D \times (0,T)$, $\Gamma_T = \{(x,y) \in \partial D, \ 0 < t < T\}$, $\partial D$ is the boundary of $D$, $D_0 = \{(x,y) \in D, t = 0\}$, $D_T = \{(x,y) \in D, t = T\}$. Then a function $u(x,y,t) \in C^{2,1}(P_T) \cap C(P_T \cup \Gamma_T \cup \bar{D}_0) \cup \Gamma_T \cup \bar{D}_0)$ satisfying the inclusion (5), the initial condition $u(x,y,0) = \alpha(x,y)$ on $D_0$, and the boundary conditions in $\Gamma_T$ is called the classical solution of the initial-boundary value problem (4) - (6). Here $C^{2,1}(P_T)$ is a space of continuous functions having in $P_T$ continuously partial derivatives $u_x(x,y,t), u_y(x,y,t), u_t(x,y,t), u_{xx}(x,y,t), u_{yy}(x,y,t)$. It should be noted that the definition of a solution in one sense or another (classical, generalized, weak, strong, fundamental, positive, etc.) in no way hinders the implementation of this method for the class of problems under consideration.

For the functions $g(\cdot, x, y, t)$, $(x, y, t) \in D \times [0,T]$ and the set-valued mapping $Q(\cdot, x, y, t)$ we impose the following condition (to avoid repetitions on the concepts of convex upper approximations (CUA) for the convenience of the reader, we refer to [25, p.122] It should only be noted that, since a nonsmooth and nonconvex function cannot be approximated in a neighbourhood of some point with positively homogeneous functions, we use the concept of CUAs.

**CONDITION N** *Let in problem (PD) the mapping $Q(\cdot, x, y, t)$ be such that the cone of tangent directions $K_{Q(\cdot,x,y,t)}(u_{x-1,y,t}, u_{x,y-1,t}, u_{x,y,t}, u_{x,y+1,t}, u_{x,y,t+1}, u_{x+1,y,t})$ is a local tent. In addition, let the functions $g(\cdot, x, y, t)$ admit CUA $h_{x,y,t}(\bar{u}, u_{x,y,t})$ at points $u_{x,y,t}$ that are continuous with respect to $\bar{u}$. The latter means that the subdifferentials $\partial g(u_{x,y,t}, x, y, t) := \partial h_{x,y,t}(0, u_{x,y,t})$ are defined. The problem (1) - (3) is called convex, if the mapping $Q(\cdot, x, y, t)$ is convex and the functions $g(\cdot, x, y, t)$, $(x, y, t) \in L_0 \times S_0 \times T_0, (x, y) \neq (0,0), (L,S)$ are proper convex functions.*

The following condition guarantees for the problem (PD) the existence of a standard interior point of convex analysis. As usual, ri $A$ be the relative interior of a set $A$.

**DEFINITION 2.1** *The convex problem (PD) satisfies the nondegenacy condition if for some feasible solution $\{u_{x,y,t}\}$ we have either (a) or (b):*

*(a)* $(u_{x-1,y,t}, u_{x,y-1,t}, u_{x,y,t}, u_{x,y+1,t}, u_{x,y,t+1}, u_{x+1,y,t}) \in \text{ri gph } Q(\cdot, x, y, t), (x, y, t) \in L_1 \times S_1 \times T_1,$

  $u_{x,y,t} \in \text{ri dom } g(\cdot, x, y, t),$

*(b)* $(u_{x-1,y,t}, u_{x,y-1,t}, u_{x,y,t}, u_{x,y+1,t}, u_{x,y,t+1}, u_{x+1,y,t}) \in \text{int gph } Q(\cdot, x, y, t), (x, y, t) \in L_1 \times S_1 \times T_1,$

$(x, y, t) \neq (x_0, y_0, t_0)$, where $(x_0, y_0, t_0)$ is the fixed triple and $g(\cdot, x, y, t)$ are continuous at the point $u_{x,y,t}$.

## 3. OPTIMIZATION OF FIRST MIXED PARABOLIC TYPE DISCRETE INCLUSIONS

We first consider the convex problem (PD). In order to use convex programming results, we form the $m = n(L+1)(S+1)(T+1)$ dimensional vector $w = (u_0, u_1, ...., u_L) \in \mathbb{R}^m$, where $u_x = (u_{x0}, u_1, ...., u_{xS}) \in \mathbb{R}^{n(S+1)(T+1)}$ and $u_{xy} = (u_{x,y,0}, u_{x,y,1}, ...., u_{x,y,T}) \in \mathbb{R}^{n(T+1)}$. Let us consider the following convex sets defined in the space $\mathbb{R}^m$:

$$M_{x,y,t} = \left\{ w = (u_0, u_1, ...., u_L) : \left( u_{x-1,y,t}, u_{x,y-1,t}, u_{x,y,t}, u_{x,y+1,t}, u_{x,y,t+1}, u_{x+1,y,t} \right) \in \mathrm{gph}\, Q(\cdot, x, y, t) \right\},$$
$$(x, y, t) \in L_1 \times S_1 \times T_1,$$
$$H_\alpha = \left\{ w = (u_0, u_1, ...., u_L) : u_{x,y,0} = \alpha_{x,y}, (x, y) \in L_0 \times S_0, (x, y) \neq (0,0), (L, S) \right\},$$
$$H_{\beta 0} = \left\{ w = (u_0, u_1, ...., u_L) : u_{x,0,t} = \beta_{x0t}, (x, t) \in L_0 \times T_0, (x, t) \neq (0,0), (L, T) \right\},$$
$$H_{\beta S} = \left\{ w = (u_0, u_1, ...., u_L) : u_{x,S,t} = \beta_{xSt}, (x, t) \in L_0 \times T_0, (x, t) \neq (0,0), (L, T) \right\},$$
$$H_{\gamma 0} = \left\{ w = (u_0, u_1, ...., u_L) : u_{0,y,t} = \gamma_{0yt}, (y, t) \in S_0 \times T_0, (y, t) \neq (0,0), (S, T) \right\},$$
$$H_{\gamma L} = \left\{ w = (u_0, u_1, ...., u_L) : u_{L,y,t} = \gamma_{Lyt}, (y, t) \in S_0 \times T_0, (y, t) \neq (0,0), (S, T) \right\}.$$

Now setting

$$g(w) = \sum_{\substack{(x,y,t) \in L_0 \times S_0 \times T_0 \\ (x,t) \neq (0,0), (L,S)}} g(u_{x,y,t}, x, y, t)$$

we reduce the convex problem (PD) to the following convex minimization problem in the space $\mathbb{R}$:

minimize $g(w)$ subject to $w \in N = \left( \bigcap_{(x,y,t) \in L_1 \times S_1 \times T_1} M_{x,y,t} \right) \bigcap H_\alpha \bigcap H_{\beta 0} \bigcap H_{\beta S} \bigcap H_{\gamma 0} \bigcap H_{\gamma L}.$ (7)

We apply Theorem 3.4 [25, p.99] to the convex minimization problem (7). For this, it is necessary to calculate the dual cones $K^*_{M_{x,y,t}}(w)$, $K^*_{H_\alpha}(w)$, $K^*_{H_{\beta 0}}(w)$, $K^*_{H_{\beta S}}(w)$, $K^*_{H_{\gamma 0}}(w)$, $K^*_{H_{\gamma L}}(w)$, $w \in N$.

**LEMMA 3.1** *The dual cone $K^*_{M_{x,y,t}}(w)$ to the cone of tangent directions $K_{M_{x,y,t}}(w)$ has a form:*

$$K^*_{M_{x,y,t}}(w) = \left\{ w^* = (u_0^*, u_1^*, ...., u_L^*) : \left( u_{x-1,y,t}^*, u_{x,y-1,t}^*, u_{x,y,t}^*, u_{x,y+1,t}^*, u_{x,y,t+1}^*, u_{x+1,y,t}^* \right) \right.$$
$$\in K^*_{\mathrm{gph}Q} \left( u_{x-1,y,t}, u_{x,y-1,t}, u_{x,y,t}, u_{x,y+1,t}, u_{x,y,t+1}, u_{x+1,y,t} \right), u_{i,j,k}^* = 0, (i,j,k) \neq (x-1, y, t),$$
$$\left. (x, y-1, t), (x, y, t), (x, y+1, t), (x, y, t+1), (x+1, y, t), (x, y, t) \in L_1 \times S_1 \times T_1 \right\}$$

*Proof.* Let $\bar{w} \in K_{M_{x,y,t}}(w), w \in N$. This means that $w + \mu \bar{w} \in M_{x,y,t}$ for sufficiently small $\mu > 0$, which is the same as

$$\left( u_{x-1,y,t} + \mu \bar{u}_{x-1,y,t}, u_{x,y-1,t} + \mu \bar{u}_{x,y-1,t}, u_{x,y,t} + \mu \bar{u}_{x,y,t}, u_{x,y+1,t} + \mu \bar{u}_{x,y+1,t}, \right.$$ (8)
$$\left. u_{x,y,t+1} + \mu \bar{u}_{x,y,t+1}, u_{x+1,y,t} + \mu \bar{u}_{x+1,y,t} \right) \in \mathrm{gph}\, Q(\cdot, x, y, t).$$

On the other hand, $w^* \in K^*_{M_{x,y,t}}(w)$ is equivalent to the condition

$$\langle \bar{w}, w^* \rangle = \sum_{(x,y,t) \in L_1 \times S_1 \times T_1} \langle \bar{u}_{x,y,t}, u^*_{x,y,t} \rangle \geq 0, \ \bar{w} \in K_{M_{x,y,t}}(w), \ w \in N,$$

where the components $\bar{u}_{x,y,t}$ of thee vector $\bar{w}$ (see(8)) are arbitrary. Therefore, the last relation is valid only for $u^*_{i,j,k} = 0$, $(i,j,k) \neq (x-1,y,t), (x,y-1,t), (x,y,t), (x,y+1,t), (x,y,t+1)$, $(x+1,y,t)$, $(x,y,t) \in L_1 \times S_1 \times T_1$. This ends the proof of the lemma. □

**THEOREM 3.1** *Assume that $Q(\cdot, x, y, t)$, $(x,y,t) \in L_1 \times S_1 \times T_1$ are convex set-valued mappings, and $g(\cdot, x, y, t)$ are convex proper functions at the points of some feasible solution $\{u_{x,y,t}\}$. Then for the $\{\tilde{u}_{x,y,t}\}$ to be an optimal solution of the problem (PD) it is necessary that there exist a number $\lambda \in \{0,1\}$ and a five of vectors $\{\psi^*_{x,y,t}, \xi^*_{x,y,t}, u^*_{x,t}, \eta^*_{x,y,t}, \varphi^*_{x,y,t}\}$, not all zero such that:*

(i) $\quad (\psi^*_{x,y,t}, \xi^*_{x,y,t}, u^*_{x-1,y,t}, \eta^*_{x,y,t}, \varphi^*_{x,y,t})$
$\in Q^*(u^*_{x,y,t}; (\tilde{u}_{x-1,y,t}, \tilde{u}_{x,y-1,t}, \tilde{u}_{x,y,t}, \tilde{u}_{x,y+1,t}, \tilde{u}_{x,y,t+1}, \tilde{u}_{x+1,y,t}), x, y, t),$
$+ \{0\} \times \{0\} \times \{\psi^*_{x+1,y,t} + \xi^*_{x,y+1,t} + \eta^*_{x,y-1,t} + \varphi^*_{x,y,t-1} - \lambda \partial g(\tilde{u}_{x,y,t}, x, y, t)\} \times \{0\} \times \{0\}$
$(x,y,t) \in L_1 \times S_1 \times T_1;$

(ii) $\quad u^*_{x,y,T} = 0, \ \eta^*_{x,0,t} = 0, \ \xi^*_{x,S,t} = 0, \ \psi^*_{1,y,t} = 0, \ \varphi^*_{L,y,t} = 0.$

*In addition, under the nondegeneracy condition, $\lambda = 1$ and these conditions are also sufficient for the optimality of $\{\tilde{u}_{x,y,t}\}$.*

*Proof.* By thehypothesis of the theorem, $\tilde{w} = (\tilde{u}_0, \tilde{u}_1, ..., \tilde{u}_L)$ is a solution of the convex minimization problem (7) and $g(w)$ is continuous at the point $w^0 = (u^0_0, ..., u^0_L)$. The we can assert the existence of vectors $w^*(x,y,t) \in K^*_{M_{x,y,t}}(\tilde{w})$, $w^*_\alpha \in K^*_{H_\alpha}(\tilde{w})$, $w^*_{\beta 0} \in K^*_{H_{\beta 0}}(\tilde{w})$, $w^*_{\beta S} \in K^*_{H_{\beta S}}(\tilde{w})$, $w^*_{\gamma 0} \in K^*_{H_{\gamma 0}}(\tilde{w})$, $w^*_{\gamma L} \in K^*_{H_{\gamma L}}(\tilde{w})$ and of a number $\lambda = 0$ or $1$, not all equal to zero, such that

$$\sum_{(x,y,t) \in L_1 \times S_1 \times T_1} w^*(x,y,t) + w^*_\alpha + w^*_{\beta 0} + w^*_{\beta S} + w^*_{\gamma 0} + w^*_{\gamma L} = \lambda w^{0*}. \tag{9}$$

On the other hand it is also easy to show that

$$K_{H_\alpha}^* = \{w^* = (u^*_0, u^*_1, ..., u^*_L): u^*_{x,y,t} = 0, \ t \neq 0, \ (x,y) \in L_0 \times S_0, (x,y) \neq (0,0), (L,S)\},$$
$$K_{H_{\beta 0}}^* = \{w^* = (u^*_0, u^*_1, ..., u^*_L): u^*_{x,y,t} = 0, \ y \neq 0, \ (x,t) \in L_0 \times T_0, (x,t) \neq (0,0), (L,T)\},$$
$$K_{H_{\beta S}}^* = \{w^* = (u^*_0, u^*_1, ..., u^*_L): u^*_{x,y,t} = 0, \ y \neq S, \ (x,t) \in L_0 \times T_0, (x,t) \neq (0,0), (L,T)\}, \tag{10}$$
$$K_{H_{\gamma 0}}^* = \{w^* = (u^*_0, u^*_1, ..., u^*_L): u^*_{x,y,t} = 0, \ x \neq 0, \ (y,t) \in S_0 \times T_0, (y,t) \neq (0,0), (S,T)\},$$
$$K_{H_{\gamma L}}^* = \{w^* = (u^*_0, u^*_1, ..., u^*_L): u^*_{x,y,t} = 0, \ x \neq L, \ (y,t) \in S_0 \times T_0, (y,t) \neq (0,0), (S,T)\}.$$

Let now $\left[ w^* \right]_{x,y,t}$ denote the components of the vector $w^*$ for the given triple $(x,y,t)$. Then using Lemma 3.1 and the relation (10), we get

$$\left[ \sum_{(x,y,t)\in L_1\times S_1\times T_1} w^*(x,y,t) + w^*_\alpha + w^*_{\beta 0} + w^*_{\beta S} + w^*_{\gamma 0} + w^*_{\gamma L} \right]_{x,y,t}$$

$$= \begin{cases} u^*_{x,y,T}(x+1,y,T) + u^{\alpha*}_{x,y,T}, & (x,y)\in L_0\times S_0, (x,y)\neq(0,0), \\ u^*_{x,S,t}(x,S-1,t) + u^{\beta S*}_{x,S,t}, & (x,t)\in L_0\times T_0, (x,t)\neq(0,0),(L,T), \\ u^*_{x,0,t}(x,1,t) + u^{\beta 0*}_{x,0,t}, & (x,t)\in L_0\times T_0, (x,t)\neq(0,0),(L,T), \\ u^*_{L,y,t}(L-1,y,t) + u^{\gamma L*}_{L,y,t}, & (y,t)\in S_0\times T_0, (y,t)\neq(0,0),(S,T), \\ u^*_{0,y,t}(1,y,t) + u^{\gamma 0*}_{0,y,t}, & (y,t)\in S_0\times T_0, (y,t)\neq(0,0),(S,T), \end{cases}$$

where it is taken into account that

$$\left[w^*_\alpha\right]_{x,y,T} = u^{\alpha*}_{x,y,T};\ \left[w^*_{\beta 0}\right]_{x,0,t} = u^{\beta 0*}_{x,0,t};\ \left[w^*_{\beta S}\right]_{x,S,t} = u^{\beta S*}_{x,S,t};\ \left[w^*_{\gamma 0}\right]_{0,y,t} = u^{\gamma 0*}_{0,y,t};\ \left[w^*_{\gamma L}\right]_{L,y,t} = u^{\gamma L*}_{L,y,t}.$$

Then because of arbitrariness of the vectors $u^{\alpha*}_{x,y,T}$, $u^{\beta_0*}_{x,0,t}$, $u^{\beta S*}_{x,S,t}$, $u^{\gamma 0*}_{0,y,t}$, $u^{\gamma L*}_{L,y,t}$ it follows from (9) that the following relationship always holds

$$u^*_{x,y,T}(x,y,T-1) + u^{\alpha*}_{x,y,T} = 0, u^*_{x,S,t}(x,S-1,t) + u^{\beta S*}_{x,S,t} = 0, u^*_{x,0,t}(x,1,t)$$
$$+ u^{\beta 0*}_{x,0,t} = 0,\ u^*_{L,y,t}(L-1,y,t) + u^{\gamma L*}_{L,y,t} = 0, u^*_{0,y,t}(1,y,t) + u^{\gamma 0*}_{0,y,t} = 0.$$

Thus, (9) implies that

$$u^*_{x,y,t}(x+1,y,t) + u^*_{x,y,t}(x,y+1,t) + u^*_{x,y,t}(x,y,t) + u^*_{x,y,t}(x,y-1,t)$$
$$+ u^*_{x,y,t}(x,y,t-1) + u^*_{x,y,t}(x-1,y,t) = \lambda u^{0*}_{x,y,t};\ \left[w^{0*}\right]_{x,y,t} = u^{0*}_{x,y,t};\ (x,y,t)\in L_1\times S_1\times T_0. \quad (11)$$
$$u^*_{x+1,y,T}(x,y,T) = 0,\ u^*_{x,S-1,t}(x,S,t) = 0,\ u^*_{x,1,t}(x,0,t) = 0,\ u^*_{L-1,y,t}(L,y,t) = 0,\ u^*_{1,y,t}(0,y,t) = 0.$$

Using Lemma 3.1 and the definition of the LAM, it can be concluded that

$$\left(u^*_{x-1,y,t}(x,y,t), u^*_{x,y-1,t}(x,y,t), u^*_{x,y,t}(x,y,t), u^*_{x,y+1,t}(x,y,t), u^*_{x,y,t+1}(x,y,t)\right)$$

$$\in Q^*\left(-u^*_{x+1,y,t}(x,y,t);\ (\tilde{u}_{x-1,y,t}, \tilde{u}_{x,y-1,t}, \tilde{u}_{x,y,t}, \tilde{u}_{x,y+1,t}, \tilde{u}_{x,y,t+1}, \tilde{u}_{x+1,y,t}), x,y,t\right),$$
$$(x,y,t)\in L_1\times S_1\times T_1\}. \quad (12)$$

Then introducing the new notations $u^*_{x-1,y,t}(x,y,t)\equiv \psi^*_{x,y,t}$, $u^*_{x,y-1,t}(x,y,t)\equiv \xi^*_{x,y,t}$, $u^*_{x,y+1,t}(x,y,t)$
$\equiv \eta^*_{x,y,t}$, $u^*_{x,y,t+1}(x,y,t)\equiv \varphi^*_{x,y,t}$, $-u^*_{x+1,y,t}(x,y,t)\equiv u^*_{x,y,t}$ we see from (11) that

$$u^*_{x,y,t}(x,y,t) = u^*_{x-1,y,t} - \psi^*_{x+1,y,t} - \xi^*_{x,y+1,t} - \eta^*_{x,y-1,t} - \varphi^*_{x,y,t-1} + \lambda u^{0*}_{x,y,t}. \quad (13)$$

Thus, taking into account these notations and substituting (13) into the Euler-Lagrange inclusion (12), we obtain the required first part of the theorem. On the other hand, it follows from nondegeneracy condition(Definition 2.1) that relation (9) holds with $\lambda = 1$ for the point $w^{0*}\in \partial g(\tilde{w}) \cap K_N^*(\tilde{w})$. Hence, condition (i) and (ii) of theorem are sufficient for the optimality of $\{\tilde{u}_{x,y,t}\}$. This completes the proof of the theorem. $\square$

**REMARK 3.1** *We note that, if in the problem (1) - (3) the functions and mappings are polyhedral, then, according to Lemma 1.22 [25, p.23] and Theorem 3.4 [25, p.99], the nondegeneracy condition in Theorem 3.1 is superfluous.*

**THEOREM 3.2** *Suppose the assumptions of condition N are carried out for the nonconvex problem (PD). Then, for the optimality of the solution $\{\tilde{u}_{x,y,t}\}$ in the nonconvex problem with a parabolic discrete inclusions, it is necessary that there exist a number $\lambda \in \{0,1\}$ and five of vectors $\{\psi^*_{x,y,t}, \xi^*_{x,y,t}, u^*_{x,t}, \eta^*_{x,y,t}, \varphi^*_{x,y,t}\}$, not all zero, satisfying the conditions (i) and (ii) of Theorem 3.1.*

*Proof.* Here the condition N ensures the conditions of Theorem 3.24 [25, p.133]. Hence, in view of this theorem, we have the necessary condition as in Theorem 3.1, written out for the nonconvex problem. □

## 4. OPTIMIZATION OF FIRST MIXED PARABOLIC TYPE DISCRETE-APPROXIMATE PROBLEM

Here, to obtain optimality conditions for the problem (PC), as an intermediate stage, we use the difference derivatives to approximate the problem (PC) and, using Theorem 3.1, we formulate a necessary and sufficient condition for it. We choose steps $\delta, \sigma$ and $h$ on the $x$, $y$-and $t$-axis, respectively, using the grid function $u_{x,y,t} = u_{\delta\sigma h}(x,y,t)$ on a uniform grid on $D \times [0,T]$. We introduce the following difference operators, defined on the three-point models [25, p.319; 26], $A_{1\delta} = A_1$ and $A_{2\sigma} = A_2$

$$A_1 u(x,y,t) = \frac{u(x+\delta,y,t) - 2u(x,y,t) + u(x-\delta,y,t)}{\delta^2};$$

$$A_2 u(x,y,t) = \frac{u(x,y+\sigma,t) - 2u(x,y,t) + u(x,y-\sigma,t)}{\sigma^2};$$

$$Bu(x,y,t) = \frac{u(x,y,t+h) - u(x,y,t)}{h}; \quad x = \delta,...,L-\delta; \ y = \sigma,...,S-\sigma; \ t = 0,h,...,T-h.$$

Therefore, we establish the following difference first mixed parabolic type boundary value problem (PDA):

$$\text{minimize} \sum_{\substack{(x,y,t) \in L_0 \times S_0 \times T_0 \\ (x,y) \neq (0,0),(L,S)}} \delta\sigma h\, g(u(x,y,t), x, y, t), \tag{14}$$

(PDA) $$Bu(x,y,t) - A_1 u(x,y,t) - A_2 u(x,y,t) \in F(u(x,y,t), x, y, t), \tag{15}$$

$$u(x,y,0) = \alpha(x,y), \ u(x,0,t) = \beta_0(x,t), \ u(x,S,t) = \beta_S(x,t),$$

$$u(0,y,t) = \gamma_0(y,t), \ u_{L,y,t} = \gamma_L(y,t), \ (x,y,t) \in \mathbb{N} \times \mathfrak{I} \times \mathfrak{R},$$

$$\mathbb{N} = \{\delta,...,L-\delta\}; \ \mathfrak{I} = \{\sigma,...,S-\sigma\}; \ \mathfrak{R} = \{0,h,...,T-h\}.$$

Recall that (14) is the Riemann sum [26] of the triple integral of a continuous function of three variables $g$ over a bounded region in three-dimensional space.

Introducing a new mapping, we reduce the problem (14) and (15) to a problem of the form (PD). To do this, rewrite discrete-approximate inclusion (15) in the equivalent form, where

$$u(x+\delta, y,t) \in -u(x-\delta, y,t) - \theta^2 u(x, y-\sigma,t) + \left(2 + 2\theta^2 - \frac{\delta^2}{h}\right) u(x, y,t)$$

$$-\theta^2 u(x, y+\sigma,t) + \frac{\delta^2}{h} u(x, y, t+h) - \delta^2 F(u(x,t), x, t), \quad \theta = \frac{\delta}{\sigma}.$$

Now let us denote

$$G\big(u(x-\delta, y,t), u(x, y-\sigma,t), u(x, y,t), u(x, y+\sigma,t), u(x, y, t+h), x, y, t\big)$$

$$= -u(x-\delta, y,t) - \theta^2 u(x, y-\sigma,t) + \left(2 + 2\theta^2 - \frac{\delta^2}{h}\right) u(x, y,t)$$

$$- \theta^2 u(x, y+\sigma,t) + \frac{\delta^2}{h} u(x, y, t+h) - \delta^2 F(u(x, y,t), x, y, t)$$

or more abbreviated

$$G(u_1, u_2, u, u_3, u_4, x, y, t) = -u_1 - \theta^2 u_2 + \left(2 + 2\theta^2 - \frac{\delta^2}{h}\right) u - \theta^2 u_3 + \frac{\delta^2}{h} u_4 - \delta^2 F(u, x, y, t) \quad (16)$$

and rewrite the problem (PDA) as follows:

$$\text{minimize} \sum_{\substack{(x,y,t) \in L_0 \times S_0 \times T_0 \\ (x,y) \neq (0,0),(L,S)}} \delta \sigma h \, g(u(x, y,t), x, y, t)$$

$$u(x+\delta, y,t) \in G\big(u(x-\delta, y,t), u(x, y-\sigma,t), u(x, y,t), u(x, y+\sigma,t), u(x, y, t+h), x, y, t\big) \quad (17)$$

$$u(x, y, 0) = \alpha(x, y), \quad u(x, 0, t) = \beta_0(x, t), \quad u(x, S, t) = \beta_S(x, t),$$

$$u(0, y, t) = \gamma_0(y, t), \quad u_{L, y, t} = \gamma_L(y, t), \quad (x, y, t) \in \mathbb{N} \times \mathfrak{I} \times \mathfrak{R}.$$

According to Theorem 3.1 for optimality of the trajectory $\{\tilde{u}(x, y, t)\}$, in problem (17) it is necessary that, there exists a five of vectors $\{\hat{\psi}^*_{x,y,t}, \hat{\xi}^*_{x,y,t}, \hat{u}^*_{x,t}, \hat{\eta}^*_{x,y,t}, \hat{\varphi}^*_{x,y,t}\}$, and a number $\lambda = \lambda_{\delta h} \in \{0, 1\}$ not all zero, such that

$$\big(\hat{\psi}^*(x, y, t), \hat{\xi}^*(x, y, t), \hat{u}^*(x-\delta, y, t), \hat{\eta}^*(x, y, t), \hat{\varphi}^*(x, y, t)\big)$$

$$\in G^*\big(\hat{u}^*(x, y, t); (\tilde{u}(x-\delta, y,t), \tilde{u}(x, y-\sigma,t), \tilde{u}(x, y,t), \tilde{u}(x, y+\sigma,t), \tilde{u}(x, y, t+h),$$

$$\tilde{u}(x+\delta, y,t)), x, y, t\big) + \{0\} \times \{0\} \times \{\hat{\psi}^*(x+\delta, y,t) + \hat{\xi}^*(x, y+\sigma,t) + \hat{\eta}^*(x, y-\sigma,t)$$

$$+ \hat{\varphi}^*(x, y, t-h) - \lambda \partial g(\tilde{u}(x, y, t), x, y, t)\} \times \{0\} \times \{0\}; \quad (x, y, t) \in \mathbb{N} \times \mathfrak{I} \times \mathfrak{R}, \quad (18)$$

$$\hat{u}^*(x, y, T) = 0, \, \hat{\eta}^*(x, 0, t) = 0, \, \hat{\xi}^*(x, S, t) = 0, \, \hat{\psi}^*(\delta, y, t) = 0, \, \hat{\varphi}^*(L, y, t) = 0. \quad (19)$$

Obviously, to form an optimality condition for the discrete-approximate problem (PDA), one must be able to transform from LAM $G^*$ to LAM $F^*$ in (18).
For this we need the following auxiliary results.

**THEOREM 4.1** Let $G(\cdot, x, y, t) : \mathbb{R}^{5n} \rightrightarrows \mathbb{R}^n$ be a set-valued mapping, defined as follows

$$G(u_1, u_2, u, u_3, u_4, x, y, t) = -u_1 - \theta^2 u_2 + \left(2 + 2\theta^2 - \frac{\delta^2}{h}\right) u - \theta^2 u_3 + \frac{\delta^2}{h} u_4 - \delta^2 F(u, x, y, t)$$

Moreover, let $G(\cdot, x, y, t)$ be a set-valued mapping such that the cone of tangent directions $K_{G(\cdot, x, y, t)}(u_1, u_2, u, u_3, u_4, v)$, $(u_1, u_2, u, u_3, u_4, v) \in \mathrm{gph} G(\cdot, x, y, t)$ is a local tent. Then

$$K_{F(\cdot, x, y, t)}\left(u, \left(\frac{2}{\delta^2} + \frac{2}{\sigma^2} - \frac{1}{h}\right)u + \frac{1}{h}u_4 - \frac{1}{\delta^2}u_1 - \frac{1}{\sigma^2}u_2 - \frac{1}{\sigma^2}u_3 - \frac{1}{\delta^2}v\right)$$

is a local tent to $\mathrm{gph} F(\cdot, x, y, t)$ and the following inclusions are equivalent:

(1)  $(\bar{u}_1, \bar{u}_2, \bar{u}, \bar{u}_3, \bar{u}_4, \bar{v}) \in K_{G(\cdot, x, y, t)}(u_1, u_2, u, u_3, u_4, v)$,

(2) $\left(\bar{u}, \left(\frac{2}{\delta^2} + \frac{2}{\sigma^2} - \frac{1}{h}\right)\bar{u} + \frac{1}{h}\bar{u}_4 - \frac{1}{\delta^2}\bar{u}_1 - \frac{1}{\sigma^2}\bar{u}_2 - \frac{1}{\sigma^2}\bar{u}_3 - \frac{1}{\delta^2}\bar{v}\right)$

$\in K_{F(\cdot, x, y, t)}\left(u, \left(\frac{2}{\delta^2} + \frac{2}{\sigma^2} - \frac{1}{h}\right)u + \frac{1}{h}u_4 - \frac{1}{\delta^2}u_1 - \frac{1}{\sigma^2}u_2 - \frac{1}{\sigma^2}u_3 - \frac{1}{\delta^2}v\right)$.

*Proof.* First let us prove the implication $(1) \Rightarrow (2)$. By the definition of a local tent there exist functions $r_i(\bar{z})$, $i = 0, 1, 2, 3, 4, 5$ and $r(\bar{z})$, $\bar{z} = (\bar{u}_1, \bar{u}_2, \bar{u}, \bar{u}_3, \bar{u}_4, \bar{v})$ such that $r_i(\bar{z})\|\bar{z}\|^{-1} \to 0$, $r(\bar{z})\|\bar{z}\|^{-1} \to 0$ as $\bar{z} \to 0$ and

$$v + \bar{v} + r_0(\bar{z}) \in -(u_1 + \bar{u}_1 + r_1(\bar{z})) - \theta^2(u_2 + \bar{u}_2 + r_2(\bar{z})) + \left(2 + 2\theta^2 - \frac{\delta^2}{h}\right)(u + \bar{u} + r(\bar{z}))$$

$$- \theta^2(u_3 + \bar{u}_3 + r_3(\bar{z})) + \frac{\delta^2}{h}(u_4 + \bar{u}_4 + r_4(\bar{z})) - \delta^2 F\left(u + \bar{u} + r(\bar{z}), x, y, t\right)$$

for sufficiently small $\bar{z} \in K$, where $K \subseteq \mathrm{ri} K_{G(\cdot, x, y, t)}(z)$ is a convex cone. Transforming this inclusion, we get

$$\left(\frac{2}{\delta^2} + \frac{2}{\sigma^2} - \frac{1}{h}\right)u + \frac{1}{h}u_4 - \frac{1}{\delta^2}u_1 - \frac{1}{\sigma^2}u_2 - \frac{1}{\sigma^2}u_3 - \frac{1}{\delta^2}v + \left(\frac{2}{\delta^2} + \frac{2}{\sigma^2} - \frac{1}{h}\right)\bar{u}$$

$$+ \frac{1}{h}\bar{u}_4 - \frac{1}{\delta^2}\bar{u}_1 - \frac{1}{\sigma^2}\bar{u}_2 - \frac{1}{\sigma^2}\bar{u}_3 - \frac{1}{\delta^2}\bar{v} + \left(\frac{2}{\delta^2} + \frac{2}{\sigma^2} - \frac{1}{h}\right)r(\bar{z})$$

$$+ \frac{1}{h}r_4(\bar{z}) - \frac{1}{\delta^2}r_1(\bar{z}) - \frac{1}{\sigma^2}r_2(\bar{z}) - \frac{1}{\sigma^2}r_3(\bar{z}) - \frac{1}{\delta^2}r_0(\bar{z}) \in F\left(u + \bar{u} + r(\bar{z}), x, y, t\right).$$

Recall that by definition of the set-valued mapping $v \in G(u_1, u_2, u, u_3, u_4, x, y, t)$, implies

$$\left(\frac{2}{\delta^2} + \frac{2}{\sigma^2} - \frac{1}{h}\right)u + \frac{1}{h}u_4 - \frac{1}{\delta^2}u_1 - \frac{1}{\sigma^2}u_2 - \frac{1}{\sigma^2}u_3 - \frac{1}{\delta^2}v \in F(u, x, y, t).$$

Therefore, from this relation it is clear that

$$\left(\bar{u}, \left(\frac{2}{\delta^2} + \frac{2}{\sigma^2} - \frac{1}{h}\right)\bar{u} + \frac{1}{h}\bar{u}_4 - \frac{1}{\delta^2}\bar{u}_1 - \frac{1}{\sigma^2}\bar{u}_2 - \frac{1}{\sigma^2}\bar{u}_3 - \frac{1}{\delta^2}\bar{v}\right)$$

$$\in K_{F(\cdot, x, y, t)}\left(u, \left(\frac{2}{\delta^2} + \frac{2}{\sigma^2} - \frac{1}{h}\right)u + \frac{1}{h}u_4 - \frac{1}{\delta^2}u_1 - \frac{1}{\delta\sigma}u_2 - \frac{1}{\delta^2}u_3 - \frac{1}{\delta^2}v\right). \qquad (20)$$

By going in the reverse direction, it is also not hard to see from (20) that $(\bar{u}_1, \bar{u}_2, \bar{u}, \bar{u}_3, \bar{u}_4, \bar{v})$
$\in K_{G(\cdot, x, y, t)}(u_1, u_2, u, u_3, u_4, v)$. This completes the proof of the theorem.  $\square$

**PROPOSITION 4.1** *Let the convex-valued mapping $G(\cdot,x,y,t)$ be given as in the formulation of Theorem 3.1. Then $G_A(u_1,u_2,u,u_3,u_4,v^*,x,y,t) = F_A(u,-v^*,x,y,t)$ and the following statements for argmaximum sets are equivalent*

(1) $v \in G_A(u_1,u_2,u,u_3,u_4,v^*,x,y,t)$;

(2) $\left(\dfrac{2}{\delta^2}+\dfrac{2}{\sigma^2}-\dfrac{1}{h}\right)u + \dfrac{1}{h}u_4 - \dfrac{1}{\delta^2}u_1 - \dfrac{1}{\sigma^2}u_2 - \dfrac{1}{\sigma^2}u_3 - \dfrac{1}{\delta^2}v \in F_A(u;-v^*,x,y,t)$.

*Proof.* First of all, we will find the connection between the Hamilton functions $H_G$ and $H_F$.

$$H_G(u_1,u_2,u,u_3,u_4,v^*) = \sup\{\langle v,v^*\rangle : v \in G(u_1,u_2,u,u_3,u_4,x,y,t)\}$$

$$= \left\langle -u_1 - \theta^2 u_2 + \left(2+2\theta^2 - \dfrac{\delta^2}{h}\right)u - \theta^2 u_3 + \dfrac{\delta^2}{h}u_4, v^*\right\rangle + \delta^2 \sup\{\langle v_1,-v^*\rangle : v_1 \in F(u,x,y,t)\}$$

$$= \left\langle -u_1 - \theta^2 u_2 + \left(2+2\theta^2 - \dfrac{\delta^2}{h}\right)u - \theta^2 u_3 + \dfrac{\delta^2}{h}u_4, v^*\right\rangle + \delta^2 H_F(u,-v^*).$$

Therefore, keeping in mind this formula for $H_G$ and $H_F$ we get the proof of proposition immediately as follows

$$G_A(u_1,u_2,u,u_3,u_4;v^*,x,y,t) = \{\, v \in G(u_1,u_2,u,u_3,u_4,x,y,t):$$

$$\langle v,v^*\rangle = H_G(u_1,u_2,u,u_3,u_4,v^*)\} = \left\{\left(\dfrac{2}{\delta^2}+\dfrac{2}{\sigma^2}-\dfrac{1}{h}\right)u + \dfrac{1}{h}u_4 - \dfrac{1}{\delta^2}u_1 - \dfrac{1}{\sigma^2}u_2 - \dfrac{1}{\sigma^2}u_3\right.$$

$$= \left\{\left(\dfrac{2}{\delta^2}+\dfrac{2}{\sigma^2}-\dfrac{1}{h}\right)u + \dfrac{1}{h}u_4 - \dfrac{1}{\delta^2}u_1 - \dfrac{1}{\sigma^2}u_2 - \dfrac{1}{\sigma^2}u_3 - \dfrac{1}{\delta^2}v \in F(u,x,t):\right.$$

$$\left\langle \left(\dfrac{2}{\delta^2}+\dfrac{2}{\sigma^2}-\dfrac{1}{h}\right)u + \dfrac{1}{h}u_4 - \dfrac{1}{\delta^2}u_1 - \dfrac{1}{\sigma^2}u_2 - \dfrac{1}{\sigma^2}u_3 - \dfrac{1}{\delta^2}v, -v^*\right\rangle$$

$$= H_F(u,-v^*)\Big\} = F_A(u;-v^*,x,y,t). \text{ The needed proof is ended.} \qquad \square$$

**THEOREM 4.2** *Assume that the mapping $G(u_1,u_2,u,u_3,u_4,x,y,t)$ is such that the cones of tangent directions $K_{G(\cdot,x,y,t)}(u_1,u_2,u,u_3,u_4,v)$ determine a local tent. Then the following inclusions are equivalent under the conditions that $u_1^* = -v^*$, $u_2^* = u_3^* = -\theta^2 v^*$, $u_4^* = (\delta^2/h)v^*$:*

(1) $(u_1^*,u_2^*,u^*,u_3^*,u_4^*) \in G^*(v^*;(u_1,u_2,u,u_3,u_4),x,y,t)$, $v \in G_A(u_1,u_2,u,u_3,u_4,v^*,x,y,t)$,

(2) $\dfrac{u^*}{\delta^2} - \left(\dfrac{2}{\delta^2}+\dfrac{2}{\sigma^2}-\dfrac{1}{h}\right)v^*$

$$\in F^*\left(-v^*;\left(u,\left(\dfrac{2}{\delta^2}+\dfrac{2}{\sigma^2}-\dfrac{1}{h}\right)u + \dfrac{1}{h}u_4 - \dfrac{1}{\delta^2}u_1 - \dfrac{1}{\sigma^2}u_2 - \dfrac{1}{\sigma^2}u_3 - \dfrac{1}{\delta^2}v\right),x,y,t\right),$$

$\left(\dfrac{2}{\delta^2}+\dfrac{2}{\sigma^2}-\dfrac{1}{h}\right)u + \dfrac{1}{h}u_4 - \dfrac{1}{\delta^2}u_1 - \dfrac{1}{\sigma^2}u_2 - \dfrac{1}{\sigma^2}u_3 - \dfrac{1}{\delta^2}v \in F_A(u;-v^*,x,y,t)$, $v^* \in \mathbb{R}^n$.

*Proof.* Suppose that the condition (1) of theorem is fulfilled. By the definition of LAM, this means that in the case of (1) we have

$$\langle \bar{u}_1,u_1^*\rangle + \langle \bar{u}_2,u_2^*\rangle + \langle \bar{u},u^*\rangle + \langle \bar{u}_3,u_3^*\rangle + \langle \bar{u}_4,u_4^*\rangle - \langle \bar{v},v^*\rangle \geq 0, \qquad (21)$$

$$(\bar{u}_1, \bar{u}_2, \bar{u}, \bar{u}_3, \bar{u}_4, \bar{v}) \in K_{G(\cdot, x, y, t)}(u_1, u_2, u, u_3, u_4, v).$$

Using the equivalence of the inclusions (1) and (2) for cones $K_{G(\cdot, x, y, t)}$ and $K_{F(\cdot, x, y, t)}$, we rewrite this inequality in a convenient form as follows

$$\langle \bar{u}, \omega \rangle - \left\langle \left( \frac{2}{\delta^2} + \frac{2}{\sigma^2} - \frac{1}{h} \right) \bar{u} + \frac{1}{h} \bar{u}_4 - \frac{1}{\delta^2} \bar{u}_1 - \frac{1}{\sigma^2} \bar{u}_2 - \frac{1}{\sigma^2} \bar{u}_3 - \frac{1}{\delta^2} \bar{v}, \psi \right\rangle \geq 0, \quad (22)$$

$$\left( \bar{u}, \left( \frac{2}{\delta^2} + \frac{2}{\sigma^2} - \frac{1}{h} \right) \bar{u} + \frac{1}{h} \bar{u}_4 - \frac{1}{\delta^2} \bar{u}_1 - \frac{1}{\sigma^2} \bar{u}_2 - \frac{1}{\sigma^2} \bar{u}_3 - \frac{1}{\delta^2} \bar{v} \right)$$

$$\in K_{F(\cdot, x, y, t)} \left( u, \left( \frac{2}{\delta^2} + \frac{2}{\sigma^2} - \frac{1}{h} \right) u + \frac{1}{h} u_4 - \frac{1}{\delta^2} u_1 - \frac{1}{\sigma^2} u_2 - \frac{1}{\sigma^2} u_3 - \frac{1}{\delta^2} v \right),$$

where $\omega$ and $\psi$ are to be determined. Carrying out the necessary transformations in (22), we get

$$\langle \bar{u}_1, \psi \rangle + \langle \bar{u}_2, \theta^2 \psi \rangle + \left\langle \bar{u}, \delta^2 \omega - \left( 2 + 2\theta^2 - \frac{\delta^2}{h} \right) \psi \right\rangle \quad (23)$$

$$+ \langle \bar{u}_3, \theta^2 \psi \rangle + \left\langle \bar{u}_4, -\frac{\delta^2}{h} \psi \right\rangle - \langle \bar{v}, -\psi \rangle \geq 0.$$

Then comparing (21), (23) it is not hard to see that we obtain the following formulas for $\omega$ and $\psi$:

$$\omega = \frac{u^*}{\delta^2} - \left( \frac{2}{\delta^2} + \frac{2}{\sigma^2} - \frac{1}{h} \right) v^*, \quad u_1^* = -v^*, \quad u_2^* = u_3^* = -\theta^2 v^*, \quad u_4^* = \frac{\delta^2}{h} v^*, \quad \psi = -v^*.$$

Thus (22) means that

$$\frac{u^*}{\delta^2} - \left( \frac{2}{\delta^2} + \frac{2}{\sigma^2} - \frac{1}{h} \right) v^*$$

$$\in F^* \left( -v^*; \left( u, \left( \frac{2}{\delta^2} + \frac{2}{\sigma^2} - \frac{1}{h} \right) u + \frac{1}{h} u_4 - \frac{1}{\delta^2} u_1 - \frac{1}{\sigma^2} u_2 - \frac{1}{\sigma^2} u_3 - \frac{1}{\delta^2} v \right), x, y, t \right).$$

Then from Theorem 3.1 we obtain condition (2), i.e., (1) $\Rightarrow$ (2). By analogy, it can be shown that (2) $\Rightarrow$ (1). As for the conditions concerning the argmaximum sets $G_A$ and $F_A$, then by Theorem 2.1 [25, p.62], they simply guarantee the non-emptiness of the LAMs $G^*$ and $F^*$. The proof of the theorem is completed. □

**REMARK 4.1** *It should be recalled that for the transition from a problem with a discrete inclusion-approximation to a continuous one, it is very important to express $G^*$ in terms of $F^*$, without which it would hardly be possible to obtain the desired result for parabolic (PC) problem. Therefore, we note that we can calculate $\partial H_{G(\cdot, x, y, t)}(u_1, u_2, u, u_3, u_4, v^*)$ and express it in terms of $\partial H_{F(\cdot, x, y, t)}(u, v^*)$ and thus using Theorem 2.1 [25, p.62], calculate $F^*$.* □

So, let us return to the Euler-Lagrange-type condition (18) and transform it into a more convenient form in what follows. Using conditions, $u_1^* = -v^*$, $u_2^* = u_3^* = -\theta^2 v^*$, $u_4^* = (\delta^2/h)v^*$, the condition (1) of Theorem 3.2, and the Euler - Lagrange type inclusion (18), it is easy to see that $\hat{\psi}^*(x,y,t) = -\hat{u}^*(x,y,t)$, $\hat{\xi}^*(x,y,t) = \hat{\eta}^*(x,y,t) = -\theta^2 \hat{u}^*(x,y,t)$, $\hat{\varphi}^*(x,y,t) = (\delta^2/h)\hat{u}^*(x,y,t)$.

Therefore, introducing these expressions for $\hat{\psi}^*(x,y,t)$, $\hat{\xi}^*(x,y,t)$, $\hat{\eta}^*(x,y,t)$, $\hat{\varphi}^*(x,y,t)$ into the Euler-Lagrange type inclusion (18) we have

$$\left( -\hat{u}^*(x,y,t), -\theta^2 \hat{u}^*(x,y,t), \hat{u}^*(x-\delta,y,t) + \hat{u}^*(x+\delta,y,t) + \theta^2 \hat{u}^*(x,y+\sigma,t) \right.$$
$$\left. + \theta^2 \hat{u}^*(x,y-\sigma,t) - \frac{\delta^2}{h}\hat{u}^*(x,y,t-h), -\theta^2 \hat{u}^*(x,y,t), \frac{\delta^2}{h}\hat{u}^*(x,y,t) \right) \quad (24)$$

$$\in G^*\left(\hat{u}^*(x,y,t); (\tilde{u}(x-\delta,y,t), \tilde{u}(x,y-\sigma,t), \tilde{u}(x,y,t), \tilde{u}(x,y+\sigma,t), \tilde{u}(x,y,t+h), \right.$$
$$\left. \tilde{u}(x+\delta,y,t)), x,y,t\right) + \{0\} \times \{0\} \times \{-\lambda \delta h \partial g(\tilde{u}(x,y,t),x,t)\} \times \{0\} \times \{0\}, \ (x,y,t) \in \mathbb{N} \times \mathfrak{I} \times \mathfrak{R}.$$

In order to pass from $G^*$ to $F^*$ here, we apply the second condition (2) of Theorem 4.2. By this condition third term in the left hand side of (24) using the difference operators $A_1, A_2$ and $B$ is converted like this:

$$\frac{1}{\delta^2}\left[ \hat{u}^*(x-\delta,y,t) + \hat{u}^*(x+\delta,y,t) + \theta^2 \hat{u}^*(x,y+\sigma,t) + \theta^2 \hat{u}^*(x,y-\sigma,t) - \frac{\delta^2}{h}\hat{u}^*(x,y,t-h) \right] \quad (25)$$

$$-\left(\frac{2}{\delta^2} + \frac{2}{\sigma^2} - \frac{1}{h}\right)\hat{u}^*(x,y,t) = \frac{\hat{u}^*(x+\delta,y,t) - 2\hat{u}^*(x,y,t) + \hat{u}^*(x-\delta,y,t)}{\delta^2}$$
$$+ \frac{\hat{u}^*(x,y+\sigma,t) - 2\hat{u}^*(x,y,t) - \hat{u}^*(x,y-\sigma,t)}{\sigma^2} + \frac{\hat{u}^*(x,y,t) - \hat{u}^*(x,y,t-h)}{h}$$
$$= A_1 \hat{u}^*(x,y,t) + A_2 \hat{u}^*(x,y,t) + B\hat{u}^*(x,y,t-h),$$

where we took into account that $\theta = \delta/\sigma$. Finally, based on (25) for the difference boundary value problem (PDA) the required Euler-Lagrange inclusion consists of the following

$$A_1 \hat{u}^*(x,y,t) + A_2 \hat{u}^*(x,y,t) + B\hat{u}^*(x,y,t-h)$$
$$\in F^*\left(-\hat{u}^*(x,t); (\tilde{u}(x,t), B\tilde{u}(x,y,t) - A_1\tilde{u}(x,y,t) - A_2\tilde{u}(x,y,t)), x,y,t\right) \quad (26)$$
$$-\lambda \delta \sigma h \partial g(\tilde{u}(x,y,t),x,y,t).$$

Recall that the LAM $F^*$ is a positive homogeneous mapping with respect to the first argument. Thus, dividing the left and right sides of (26) by $\delta \sigma h$ and then denoting $-\hat{u}^*(x,y,t)/\delta\sigma h \equiv u^*(x,t)$, we finally obtain for the problem (PDA) the following Euler-Lagrange type adjoint inclusion

$$-A_1 u^*(x,y,t) - A_2 u^*(x,y,t) - Bu^*(x,y,t-h)$$
$$\in F^*\left(u^*(x,t); (\tilde{u}(x,t), B\tilde{u}(x,y,t) - A_1\tilde{u}(x,y,t) - A_2\tilde{u}(x,y,t)), x,y,t\right) \quad (27)$$
$$-\lambda \partial g(\tilde{u}(x,y,t),x,y,t), \ (x,y,t) \in \mathbb{N} \times \mathfrak{I} \times \mathfrak{R},$$

where

$$u^*(x,y,T) = 0, \ u^*(x,0,t) = 0, \ u^*(x,S,t) = 0, \ u^*(\delta,y,t) = 0, \ u^*(L,y,t) = 0. \quad (28)$$

Let us state the result obtained.

**THEOREM 4.3** *Let $F(\cdot,x,y,t)$ be a convex set-valued mapping, and $g(\cdot,x,y,t)$ be a proper*

*convex function. Then for the optimality of the solution* $\{\tilde{u}(x,y,t)\}$ *in the discrete-approximate problem (PDA) it is necessary that there exist a number* $\lambda \in \{0,1\}$ *and vectors* $\{u^*(x,y,t)\}$, *not all zero, satisfying (27) and (28). In addition, under conditions of nondegeneracy, (27) and (28) are also sufficient for the optimality of* $\{\tilde{u}(x,y,t)\}$.

**REMARK 4.2** *If we consider the nonconvex problem (1)-(3) for a first mixed parabolic discrete inclusions, then it is easy to see that, under the condition N, the Euler-Lagrange type inclusion (27) and the initial-boundary values (28) may not be sufficient optimality conditions. This is explained by the fact that when forming the optimality condition in the nonconvex case, instead of Theorem 3.2[25, p.98], we use Theorems 3.24 [25, p.133], where the equality* $\lambda = 1$ *is optional.*

## 5. SUFFICIENT CONDITIONS OF OPTIMALITY FOR THE PARABOLIC DFIS

In this section, relying on the results obtained in Section 4, we formulate a sufficient optimality condition for the continuous problem (PC). To do this, recall that the limits of $A_1 v(x,y,t), A_2 v(x,y,t)$ and $B v(x,y,t)$, are first and second order partial derivatives with respect to $x, y$ and $t$, respectively, i.e., $\lim_{\delta \to 0} A_1 v(x,y,t) = \dfrac{\partial^2 v(x,y,t)}{\partial x^2}$, $\lim_{\sigma \to 0} A_2 v(x,y,t) = \dfrac{\partial^2 v(x,y,t)}{\partial y^2}$ and $\lim_{h \to 0} B v(x,y,t) = \dfrac{\partial v(x,y,t)}{\partial t}$. Thus, setting $\lambda = 1$ in (19), (20) and passing to the formal limit, as the discrete steps $\delta, \sigma$ and $h$ tend to 0, in term of Laplace operator we find that:

(i) $-\dfrac{\partial u^*(x,y,t)}{\partial t} - \Delta u^*(x,y,t) \in F^*\left(u^*(x,y,t); \left(\tilde{u}(x,y,t), \dfrac{\partial \tilde{u}(x,y,t)}{\partial t} - \Delta \tilde{u}(x,y,t)\right), x,y,t\right)$

$-\partial g\left(\tilde{u}(x,y,t),x,y,t\right), \quad (x,y,t) \in D \times [0,T]$,

(ii) $u^*(x,y,T) = u^*(x,0,t) = u^*(x,S,t) = u^*(0,y,t) = u^*(L,y,t) = 0$.

In fact, the following condition ensures that the LAM $F^*$ is nonempty (Theorem 2.1 [25, p.62] that:

(iii) $\dfrac{\partial \tilde{u}(x,y,t)}{\partial t} - \Delta \tilde{u}(x,y,t) \in F_A\left(\tilde{u}(x,y,t); u^*(x,y,t), x, y, t\right)$.

Here we assume that feasible solutions $u^*(x,y,t)$ are classical solutions satisfying an adjoint parabolic inclusion of the Euler-Lagrange type (*i*) and homogeneous end-point and boundary conditions (*ii*).

As expected, in the next theorem we prove that conditions (i) - (iii) are sufficient conditions for optimality in the first mixed problem with parabolic DFIs (PC) in a three space dimension.

**THEOREM 5.1:** *Let* $g(\cdot, x, y, t)$ *be continuous convex proper function and* $F(\cdot, x, y, t): \mathbb{R}^n \rightrightarrows \mathbb{R}^n$ *be a convex set-valued mapping. Then for the optimality of the solution* $\tilde{u}(x,y,t)$ *in the first mixed problem (PC) with parabolic DFIs and Laplace operator it is sufficient that there exists a classical solutions* $u^*(x,y,t)$, *satisfying the conditions (i)-(iii).*

*Proof.* By virtue of the definition of LAM and the definition of the Hamiltonian function $H_F$, the Euler-Lagrange type (*i*) gains the form:

$$H_F\left(u(x,y,t);u^*(x,y,t)\right) - H_F\left(\tilde{u}(x,y,t);u^*(x,y,t)\right)$$

$$\leq \left\langle -\frac{\partial u^*(x,y,t)}{\partial t} - \Delta u^*(x,y,t), u(x,y,t) - \tilde{u}(x,y,t) \right\rangle$$

$$+ g\left(u(x,y,t),x,y,t\right) - g\left(\tilde{u}(x,y,t),x,y,t\right)$$

from which under the condition (*iii*) we get

$$\left\langle \frac{\partial u(x,y,t)}{\partial t} - \Delta u(x,y,t), u^*(x,y,t) \right\rangle - \left\langle \frac{\partial \tilde{u}(x,y,t)}{\partial t} - \Delta \tilde{u}(x,y,t), u^*(x,y,t) \right\rangle$$

$$\leq \left\langle -\frac{\partial u^*(x,y,t)}{\partial t} - \Delta u^*(x,y,t), u(x,y,t) - \tilde{u}(x,y,t) \right\rangle$$

$$+ g\left(u(x,y,t),x,y,t\right) - g\left(\tilde{u}(x,y,t),x,y,t\right)$$

or

$$g\left(u(x,y,t),x,y,t\right) - g\left(\tilde{u}(x,y,t),x,y,t\right)$$

$$\geq \left\langle \frac{\partial\left(u(x,y,t) - \tilde{u}(x,y,t)\right)}{\partial t} - \Delta\left(u(x,y,t) - \tilde{u}(x,y,t)\right), u^*(x,y,t) \right\rangle$$

$$+ \left\langle \frac{\partial u^*(x,y,t)}{\partial t} + \Delta u^*(x,y,t), u(x,t) - \tilde{u}(x,t) \right\rangle.$$

We rewrite this inequality in the following form

$$g\left(u(x,y,t),x,y,t\right) - g\left(\tilde{u}(x,y,t),x,y,t\right)$$

$$\geq \left\langle \frac{\partial\left(u(x,y,t) - \tilde{u}(x,y,t)\right)}{\partial t}, u^*(x,y,t) \right\rangle + \left\langle \frac{\partial u^*(x,y,t)}{\partial t}, u(x,t) - \tilde{u}(x,t) \right\rangle$$

$$+ \left\langle \Delta u^*(x,y,t), u(x,t) - \tilde{u}(x,t) \right\rangle - \left\langle \Delta\left(u(x,y,t) - \tilde{u}(x,y,t)\right), u^*(x,y,t) \right\rangle. \tag{29}$$

Notice that it is easy to check the following formulas

$$\frac{\partial}{\partial t}\left\langle u(x,y,t) - \tilde{u}(x,y,t), u^*(x,y,t) \right\rangle = \left\langle \frac{\partial\left(u(x,y,t) - \tilde{u}(x,y,t)\right)}{\partial t}, u^*(x,y,t) \right\rangle \tag{30}$$

$$+ \left\langle u(x,y,t) - \tilde{u}(x,y,t), \frac{\partial u^*(x,y,t)}{\partial t} \right\rangle,$$

$$\left\langle \Delta u^*(x,y,t), u(x,t) - \tilde{u}(x,t) \right\rangle - \left\langle \Delta\left(u(x,y,t) - \tilde{u}(x,y,t)\right), u^*(x,y,t) \right\rangle$$

$$= \left\langle \frac{\partial^2 u^*(x,y,t)}{\partial x^2} + \frac{\partial^2 u^*(x,y,t)}{\partial y^2}, u(x,y,t) - \tilde{u}(x,y,t) \right\rangle$$

$$- \left\langle \frac{\partial^2\left(u(x,y,t) - \tilde{u}(x,y,t)\right)}{\partial x^2} + \frac{\partial^2\left(u(x,y,t) - \tilde{u}(x,y,t)\right)}{\partial y^2}, u^*(x,y,t) \right\rangle$$

$$= \left\langle \frac{\partial^2 u^*(x,y,t)}{\partial x^2}, u(x,y,t) - \tilde{u}(x,y,t) \right\rangle - \left\langle \frac{\partial^2\left(u(x,y,t) - \tilde{u}(x,y,t)\right)}{\partial x^2}, u^*(x,y,t) \right\rangle \tag{31}$$

$$= \frac{\partial}{\partial x}\left[\left\langle u(x,y,t)-\tilde{u}(x,y,t), \frac{\partial u^*(x,y,t)}{\partial x}\right\rangle - \left\langle \frac{\partial(u(x,y,t)-\tilde{u}(x,y,t))}{\partial x}, u^*(x,y,t)\right\rangle\right]$$

$$+ \frac{\partial}{\partial y}\left[\left\langle u(x,y,t)-\tilde{u}(x,y,t), \frac{\partial u^*(x,y,t)}{\partial y}\right\rangle - \left\langle \frac{\partial(u(x,y,t)-\tilde{u}(x,y,t))}{\partial y}, u^*(x,y,t)\right\rangle\right].$$

Then in view of (30) and (31) integrating the inequality (29) over the region $D \times [0,T]$, we have

$$\int_0^T \iint_D \left[g(u(x,y,t),x,y,t) - g(\tilde{u}(x,y,t),x,y,t)\right] dxdydt$$

$$\geq \int_0^T \iint_D \left[\frac{\partial}{\partial t}\left\langle u(x,y,t)-\tilde{u}(x,y,t), u^*(x,y,t)\right\rangle\right] dxdydt \qquad (32)$$

$$+ \int_0^T \iint_D \frac{\partial}{\partial x}\left[\left\langle u(x,y,t)-\tilde{u}(x,y,t), \frac{\partial u^*(x,y,t)}{\partial x}\right\rangle - \left\langle \frac{\partial(u(x,y,t)-\tilde{u}(x,y,t))}{\partial x}, u^*(x,y,t)\right\rangle\right] dxdydt$$

$$+ \int_0^T \iint_D \frac{\partial}{\partial y}\left[\left\langle u(x,y,t)-\tilde{u}(x,y,t), \frac{\partial u^*(x,y,t)}{\partial y}\right\rangle - \left\langle \frac{\partial(u(x,y,t)-\tilde{u}(x,y,t))}{\partial y}, u^*(x,y,t)\right\rangle\right] dxdydt.$$

For brevity of notation, denoting first, second and third integrals on the right hand-side of (32) by $J_1$, $J_2$ and $J_3$, respectively, we deduce

$$J_1 = \iint_D \left[\left\langle u(x,y,T)-\tilde{u}(x,y,T), u^*(x,y,T)\right\rangle\right] dxdy$$

$$- \iint_D \left[\left\langle u(x,y,0)-\tilde{u}(x,y,0), u^*(x,y,0)\right\rangle\right] dxdy.$$

Here since $u(x,y,t)$, $\tilde{u}(x,y,t)$ are feasible, by the conditions $u(x,y,0) = \tilde{u}(x,y,0) = \alpha(x,y)$ and $u^*(x,y,T) = 0$ (see (ii)) it follows that $J_1 = 0$.

Let us evaluate the integral $J_2$.

$$J_2 = \int_0^T \iint_D \frac{\partial}{\partial x}\left[\left\langle u(x,y,t)-\tilde{u}(x,y,t), \frac{\partial u^*(x,y,t)}{\partial x}\right\rangle\right.$$

$$\left. - \left\langle \frac{\partial(u(x,y,t)-\tilde{u}(x,y,t))}{\partial x}, u^*(x,y,t)\right\rangle\right] dxdydt \qquad (33)$$

$$= \int_0^S \int_0^T \left[\left\langle u(L,y,t)-\tilde{u}(L,y,t), \frac{\partial u^*(L,y,t)}{\partial x}\right\rangle - \left\langle \frac{\partial(u(L,y,t)-\tilde{u}(L,y,t))}{\partial x}, u^*(L,y,t)\right\rangle\right] dydt$$

$$- \int_0^S \int_0^T \left[\left\langle u(0,y,t)-\tilde{u}(0,y,t), \frac{\partial u^*(0,y,t)}{\partial x}\right\rangle - \left\langle \frac{\partial(u(0,y,t)-\tilde{u}(0,y,t))}{\partial x}, u^*(0,y,t)\right\rangle\right] dydt.$$

Hence, by analogy, since $u(0,y,t) = \tilde{u}(0,y,t) = \beta_0(x,t)$, $u(L,y,t) = \tilde{u}(L,y,t) = \gamma_L(y,t)$ and by condition (ii) $u^*(0,y,t) = u^*(L,y,t) = 0$, from equality (33) we obtain $J_2 = 0$.

It can be shown in exactly the same way that $J_3 = 0$. Indeed

$$J_3 = \int_0^L \int_0^T \left[ \left\langle u(x,S,t) - \tilde{u}(x,S,t), \frac{\partial u^*(x,S,t)}{\partial y} \right\rangle - \left\langle \frac{\partial \left(u(x,S,t) - \tilde{u}(x,S,t)\right)}{\partial y}, u^*(x,S,t) \right\rangle \right] dxdt$$

$$- \int_0^L \int_0^T \left[ \left\langle u(x,0,t) - \tilde{u}(x,0,t), \frac{\partial u^*(x,0,t)}{\partial y} \right\rangle - \left\langle \frac{\partial \left(u(x,0,t) - \tilde{u}(x,0,t)\right)}{\partial y}, u^*(x,0,t) \right\rangle \right] dxdt. \quad (34)$$

Therefore, since $u(x,0,t) = \tilde{u}(x,0,t) = \beta_0(x,t)$, $u(x,S,t) = \tilde{u}(x,S,t) = \beta_S(x,t)$ and $u^*(x,0,t) = u^*(x,S,t) = 0$ from equality (34) we obtain $J_3 = 0$.

Finally, from inequality (32) we have $\int_0^T \iint_D \left[ g\left(u(x,y,t), x, y, t\right) - g\left(\tilde{u}(x,y,t), x, y, t\right) \right] dxdydt \geq 0$ or $J[u(\cdot,\cdot,\cdot)] \geq J[\tilde{u}(\cdot,\cdot,\cdot)]$ for all feasible $u(\cdot,\cdot,\cdot)$. The theorem is proved. □

**REMARK 5.1:** *It is well known [30] that in the problem with ordinary polynomial linear differential operators $Ax = \sum_{k=1}^s p_k(t) D^k x$ of the $s$-th order with variable coefficients $p_k : [0,T] \to \mathbb{R}^1$ and with the operator of derivatives $D^k$ ($k = 1,...,s$) of the k-th order, the Euler-Lagrange inclusion includes the adjoint operator $A^* x^*(t) = \sum_{k=1}^s (-1)^k D^k \left[ p_k(t) x^*(t) \right]$, where the sign of the corresponding term does not change for even k. The Euler-Lagrange inclusion (i) shows that the same property holds for the second-order parabolic operator, and therefore only the partial derivative of the first order of $u^*(x,y,t)$ with respect to t changes sign.* □

Now using Theorem 5.1, consider the following linear problem as an example:

$$\text{minimize} \quad J[u(\cdot,\cdot,\cdot)] = \int_0^T \iint_D g\left(u(x,y,t), x, y, t\right) dxdydt$$

$$\frac{\partial u(x,y,t)}{\partial t} - \Delta u(x,y,t) = Au(x,y,t) + Bw(x,u,t), \quad w(x,y,t) \in U, \quad (35)$$

$$u(x,y,0) = \alpha(x,y), \quad u(x,0,t) = \beta_0(x,t), \quad u(x,S,t) = \beta_S(x,t),$$

$$u(0,y,t) = \gamma_0(y,t), \quad u(L,y,t) = \gamma_L(y,t),$$

where $A$ and $B$ are $n \times n$ and $n \times r$ matrices, respectively, $U \subset \mathbb{R}^r$ is a convex closed set, and $g$ is continuously differentiable function of $u$. It is required to find a control parameter $\tilde{w}(x,y,t) \in U$, $(x,y,t) \in D \times [0,T]$ that minimizes $J[u(\cdot,\cdot,\cdot)]$. Here $F(u) = Au + BU$. It is easy to see that

$$F^*(v^*, (u,v)) = \begin{cases} A^* v^*, & -B^* v^* \in [\text{cone}(U - w)]^*, \\ \varnothing, & -B^* v^* \notin [\text{cone}(U - w)]^*, \end{cases} \quad (36)$$

$$v = Au + Bw.$$

Therefore, using (36), and Theorem 5.1, we get

$$-\frac{\partial u^*(x,y,t)}{\partial t} - \Delta u^*(x,y,t) = A^* u^*(x,y,t) - g'_u(\tilde{u}(x,y,t), x, y, t), \quad (37)$$

$$u^*(x,y,T) = u^*(x,0,t) = u^*(x,S,t) = u^*(0,y,t) = u^*(L,y,t) = 0. \quad (38)$$

Next, we have
$$\langle w - \tilde{w}(x,y,t), B^*u^*(x,y,t)\rangle \leq 0, \quad w \in U,$$
which implies the Pontryagin maximum principle [35, p.112]:
$$\langle B\tilde{w}(x,y,t), u^*(x,y,t)\rangle = \max_{w\in U}\langle Bw, u^*(x,y,t)\rangle. \tag{39}$$
Thus, we have obtained the following result.

**THEOREM 5.2:** *The solution $\tilde{u}(x,y,t)$ of parabolic type linear problem (35) corresponding to the control function $\tilde{w}(x,y,t)$ minimizes $J[u(\cdot,\cdot,\cdot)]$, if there exists a solution $u^*(x,y,t)$ of adjoint equation (37), satisfying (38), (39).*

**COROLLARY 5.1:** *Suppose that we have an initial-boundary value problem (4) - (6), where $F(u,x,t) \equiv U$ is a constant mapping, $U \subset \mathbb{R}^n$ is a convex compact set, $g(\cdot,x,y,t)$ is a continuously differentiable function. Then the control function $\tilde{w}(x,y,t) \in U$ minimizes $J[u(\cdot,\cdot,\cdot)]$, if the solution $u^*(x,y,t)$ of the adjoint equation*
$$\frac{\partial u^*(x,y,t)}{\partial t} + \Delta u^*(x,y,t) = g'_u(\tilde{u}(x,y,t),x,y,t),$$
*with homogeneous end-point and boundary conditions (ii) of Theorem 5.1, satisfies the maximum principle*
$$\langle \tilde{w}(x,y,t), u^*(x,y,t)\rangle = \max_{w\in U}\langle w, u^*(x,y,t)\rangle.$$

*Proof.* It is obvious that $\mathrm{gph}F = \mathbb{R}^n \times U$ and so $K_F(u,v) = \mathbb{R}^n \times \mathrm{cone}(U-v)$, whence
$$K_F^*(u,v) = \{0\} \times [\mathrm{cone}(U-v)]^*.$$
Therefore, we find that
$$F^*(v^*;(u,v)) = \begin{cases} 0, & \text{if } v^* \in [\mathrm{cone}(U-v)]^*, \\ \varnothing, & \text{if } v^* \notin [\mathrm{cone}(U-v)]^*. \end{cases}$$
Then, using Theorem 5.1, we obtain the required result:
$$-\frac{\partial u^*(x,y,t)}{\partial t} - \Delta u^*(x,y,t) = -g'_u(\tilde{u}(x,y,t),x,y,t); \quad \langle w - \tilde{w}(x,y,t), u^*(x,y,t)\rangle \leq 0, \quad w \in U.$$
Note that the same result can be obtained from (37) and (39), assuming that $A$ is a zero matrix, $B$ is $n \times n$ unit matrix and $U \subset \mathbb{R}^n$. □

**REMARK 5.2:** *In particular, if $n=1$, $U = [-1,+1]$ and $g(u,x,t) \equiv u$, then by Weierstrass-Pontryagin maximum condition $\tilde{w}(x,y,t) \cdot u^*(x,y,t) = \max_{-1 \leq w \leq 1} w \cdot u^*(x,y,t)$, whence $\tilde{w}(x,y,t) = \mathrm{sgn}\, u^*(x,y,t)$, that is $\tilde{w}(x,y,t) = 1$, if $u^*(x,y,t) > 0$ and $\tilde{w}(x,y,t) = -1$, if $u^*(x,y,t) < 0$. Therefore, since $g'_u(u(x,y,t),x,y,t) = 1$ the solution $\tilde{u}(x,y,t)$ corresponding to the control function $\tilde{w}(x,y,t) = \pm 1$ minimizes $J[u(\cdot,\cdot,\cdot)]$, if $u^*(x,y,t)$ is a solution of the adjoint equation*
$$\frac{\partial u^*(x,y,t)}{\partial t} + \Delta u^*(x,y,t) = 1.$$
□

## 6. MODEL OF PARABOLIC PROBLEMS WITH POLYHEDRAL SET-VALUED MAPPING

Now consider the following problem with the so-called parabolic polyhedral DFI:

$$\text{minimize} \quad J[u(\cdot,\cdot,\cdot)] = \int_0^T \iint_D g(u(x,y,t),x,y,t)\,dxdydt,$$

(PPC)
$$A\left(\frac{\partial u(x,y,t)}{\partial t} - \Delta u(x,y,t)\right) - Bu(x,y,t) \leq d, \quad (x,y,t) \in D \times [0,T],$$

$$u(x,y,0) = \alpha(x,y), \ u(x,0,t) = \beta_0(x,t), \ u(x,S,t) = \beta_S(x,t),$$

$$u(0,y,t) = \gamma_0(y,t), \ u(L,y,t) = \gamma_L(y,t),$$

where $A, B$ are $s \times n$ dimensional matrices, $d$ is a $s$-dimensional column-vector, $g(\cdot, x, y, t)$ is a polyhedral function. Here $F(u) = \{v : Au - Bv \leq d\}$. Now we will use Theorem 5.1 for problem (PPC). By analogy formula (2.51) [25, p.84] we derive that

$$F^*(v^*;(\tilde{u},\tilde{v})) = \{-A^*q : v^* = -B^*q, q \geq 0, \langle A\tilde{u} - B\tilde{v} - d, q \rangle = 0\}.$$

Then applying Theorem 5.1 we get

$$\frac{\partial u^*(x,y,t)}{\partial t} + \Delta u^*(x,y,t) - A^*q(x,y,t) \in \partial g(\tilde{u}(x,y,t),x,y,t), \ u^*(x,y,t) = -B^*q(x,y,t),$$

$$q(x,y,t) \geq 0, \ \left\langle A\left(\frac{\partial \tilde{u}(x,y,t)}{\partial t} - \Delta \tilde{u}(x,y,t)\right) - B\tilde{u}(x,y,t) - d, q(x,y,t) \right\rangle = 0,$$

$$u^*(x,y,T) = u^*(x,0,t) = u^*(x,S,t) = u^*(0,y,t) = u^*(L,y,t) = 0. \tag{40}$$

Differentiating $u^*(x,y,t) = -B^*q(x,y,t)$ and substituting it into the first relation (40), we obtain

$$B^*\left(\Delta q(x,y,t) - \frac{\partial q(x,y,t)}{\partial t}\right) - A^*q(x,y,t) \in \partial g(\tilde{u}(x,y,t),x,y,t), \tag{41}$$

$$\left\langle A\left(\frac{\partial \tilde{u}(x,y,t)}{\partial t} - \Delta \tilde{u}(x,y,t)\right) - B\tilde{u}(x,y,t) - d, q(x,y,t) \right\rangle = 0$$

$$B^*q(x,y,T) = B^*q(x,0,t) = B^*q(x,S,t) = B^*q(0,y,t)$$
$$= B^*q(L,y,t) = 0, \ q(x,y,t) \geq 0, \ (x,y,t) \in D \times [0,T] \tag{42}$$

Thus, we obtain the following theorem

**THEOREM 6.2** *Let $g(\cdot,x,t): \mathbb{R}^n \to \mathbb{R}^1$ be a polyhedral function and that $F: \mathbb{R}^n \rightrightarrows \mathbb{R}^n$ be a polyhedral set-valued mapping. Then for the optimality of the solution $u(x,y,t), (x,y,t) \in D \times [0,T]$ in problem (PPC) with parabolic polyhedral DFIs, it is sufficient that there exists a nonnegative function $q(x,y,t) \geq 0$, $(x,y,t) \in D \times [0,T]$ satisfying the partial DFI of the Euler-Lagrange type (41) and the conjugate boundary condition (42).* □

## 7. CONCLUSION

Optimality conditions for parabolic DFIs with Laplace operator are formulated on the basis of LAM and a discretized method. For simplicity of presentation, the paper mainly considers convex problems, but with the use of local tent and CUAs, ways of extending the results obtained to the non-convex case are indicated. The main tools for constructing optimal conditions are Euler-Lagrange inclusions. In this case, an essential role is also played by equivalence theorems characteristic of problems with parabolic discrete inclusions and parabolic DFIs with Laplace operator, connecting problems (PDA) and (PC). It should also be noted that the formulation of sufficient conditions for a parabolic PDF is obtained by passing to the formal limit as the discrete steps tend to zero. Another problem of discussion, of course, is the convergence of the indicated limiting procedure and the establishment of the necessary optimality conditions, as well as the reduction of the "gap" between the necessary and sufficient conditions. In addition, the indicated difference problems, which are of independent interest, can play an important role in computational procedures.

# REFERENCES

[1] N. ABADA, M. BENCHOHRA AND H.HAMMOUCHE, *Existence and controllability results for impulsive partial functional differential inclusions*, Nonlin. Anal. Theory, Methods Appl. 69 (2008), pp. 2892-2909.

[2] J-P. AUBIN, *Viability Theory*, Birkhäuser, Springer, Boston, MA, 2009.

[3] J-P. AUBIN, *Boundary-Value Problems for Systems of Hamilton-Jacobi-Bellman Inclusions with Constraints,* SIAM J. Contr. Optim. 41 (2002), pp. 425-456

[4] H.T.BANKS AND K.KUNISCH, *An Approximation Theory for Nonlinear Partial Differential Equations with Applications to Identification and Control*, SIAM J. Contr. Optim. 20 (1982), pp.815-849.

[5] G.DI BLASIO, K.KUNISCH AND E.SINESTRARI, *L2-Regularity for Parabolic Partial Integrodifferential Equations with Delay in the Highest Order Derivatives,* J.Math. Anal. Appl. 102(1984), pp. 38-57 .

[6] P.CANNARSA AND KHAI T. NGUYEN, *Exterior Sphere Condition and Time Optimal Control for Differential Inclusions,* SIAM J. Contr. Optim. 49 (2011), pp. 2558-2576.

[7] A.CERNEA, *Some second-order necessary conditions for nonconvex hyperbolic differential inclusion problem.* J. Math. Anal. Appl. 253(2001), pp. 616-639.

[8] F. H. CLARKE, *Optimization and Nonsmooth Analysis*, New York: John Wiley, 1983.

[9] F.S.DE BLASI AND G.PIANIGIANI, *Baire category and boundary value problems for ordinary and partial differential inclusions under Caratheodory assumptions*. J. Diff. Equ. 243 (2007), pp. 558-577.

[10] A. L. DONTCHEV AND V. M. VELIOV, *Singular perturbations in Mayer's problem for linear systems.* SIAM J. Contr. Optim. 21(1983), pp. 566-581.

[11] A. L. DONTCHEV AND F. LEMPIO, *Difference methods for differential inclusions - a survey,* SIAM Rev. 34 (1992), pp. 263-294.

[12] T.DONCHEV, *Singularly Perturbed Evolution Inclusions,* SIAM J. Contr. Optim. 48 (2010), pp. 4572-4590.

[13] T.DONCHEV AND I. SLAVOV, *Averaging Method for One-Sided Lipschitz Differential Inclusions with Generalized Solutions*, SIAM J. Contr. Optim. 37(1999), pp. 1600-1613.

[14] H. FRANKOWSKA AND I. HAIDAR, *A relaxation result for state constrained delay differential inclusions*, IEEE Transactions on Automatic Control, 63(2018), pp.3751-3760.

[15] M.FRIGON, *On a critical point theory for multivalued functionals and applications to partial differential inclusions*. Nonlin. Anal. 31(1998), pp. 735-753.

[16] R. GOEBEL, *Set-valued Lyapunov functions for difference inclusions*, Automatica 47 (2011), pp.127-132,

[17] M. KISIELEWICZ, *Some optimal control problems for partial differential inclusions*. Opuscula Mathematica, 28 (2008). pp. 507-516.

[18] K.KUNISCH AND G. PEICHL, *On the Shape of Solutions of Second Order Parabolic Partial Differential Equations*, J. Diff. Equ. 75 (1988), pp.329-353.

[19] A. B. Kurzhanski and V. M. Veliov, *Set-Valued Analysis and Differential Inclusions (Progress in Systems and Control Theory, 16),* Birkhäuser; 1993.

[20] P. D. Loewen and R. T. Rockafellar, *Optimal Control of Unbounded Differential Inclusions*, SIAM J. Contr. Optim. 32 (1994), pp. 442-470.

[21] E. N. Mahmudov, *Optimization of discrete inclusions with distributed parameters*, Optimization, A Journal of Math. Program. Oper. Research, 21(1990), pp.197-207.

[22] E. N. Mahmudov, *Necessary and sufficient conditions for discrete and differential inclusions of elliptic type*. J. Math. Anal. Appl. 323 (2006), pp.768-789.

[23] E.N. Mahmudov, *Locally adjoint mappings and optimization of the first boundary value problem for hyperbolic type discrete and differential inclusions*. Nonlin. Anal. 67(2007), pp. 2966–2981.

[24] E. N. Mahmudov, *Mathematical programming and polyhedral optimization of second order discrete and differential inclusions,* Pacific J. Optim. 11(2010), pp. 511-525.

[25] E. N. Mahmudov, *Approximation and optimization of discrete and differential inclusions*. Elsevier, Boston, USA, 2011.

[26] E. N. Mahmudov, *Single variable differential and integral calculus, mathematical analysis*, Springer, Paris, France, 2013.

[27] E. N. Mahmudov, *Approximation and optimization of Darboux type differential inclusions with set-valued boundary conditions*. Optim. Letters, 7(2013), pp.871-891.

[28] E. N. Mahmudov, *Approximation and Optimization of Higher order discrete and differential inclusions*, Nonlin. Diff. Equ. Appl. 21 (2014), pp.1-26.

[29] E. N. Mahmudov, *Optimization of Mayer problem with Sturm-Liouville-type differential inclusions,* J. Optim. Theory Appl. 177 (2018), pp.345-375.

[30] E. N. Mahmudov, *Optimal control of evolution differential inclusions with polynomial linear differential operators*, Evol. Equ. Contr. Theory 8 (2019), pp.603-619.

[31] E. N. Mahmudov, *Optimal control of higher order differential inclusions with functional constraints.* ESAIM: COCV. (2019) pp.1-24. https://doi.org/ 10.1051/ cocv/ 2019018.

[32] E.N. Mahmudov, M.J. Mardanov, *On duality in optimal control problems with second-order differential inclusions and initial-point constraints*, Proc. Inst. Math. Mech. Natl. Acad. Sci. Azerb. 46 (2020), pp.115-128.

[33] B. S. Mordukhovich, *Variational Analysis and Generalized Differentiation*, I: Basic Theory; II: Applications, Grundlehren Series (Fundamental Principles of Mathematical Sciences), Vol. 330 and 331, Springer-Verlag Berlin Heidelberg, 2006.

[34] B. S. Mordukhovich, N. M. Nam, R. B. Rector and T. Tran, *Variational geometric approach to generalized differential and conjugate calculi in convex analysis*, Set-Valued Var. Anal. 25 (2017), pp. 731-755.

[35] B. S. Mordukhovich and L. Wang, *Optimal Control of Neutral Functional-Differential Inclusions*, SIAM J. Contr. Optim. 43 (2004), pp. 111-136.

[36] L. S. Pontryagin, V.G.Boltyanskii, R. V. Gamkrelidze and E. F. Mishchenko. *The mathematical theory of optimal processes*. John Wiley & Sons, Inc., New York,1962.

[37] M. Quincampoix and V. M. Veliov, *Solution tubes to differential inclusions within a collection of sets*, Contr. Cyber. 31 (2002), pp.847-862.

[38] R. T. Rockafellar and P. R. Wolenski, *Convexity in Hamilton-Jacobi theory 1: Dynamics and duality*, SIAM J. Contr. Optim., 40 (2001), pp. 1323-1350.

[39] J. Shen and Jianghi Hu, *Stability of Discrete-Time Switched Homogeneous Systems on Cones and Conewise Homogeneous Inclusions*, SIAM J. Contr. Optim. 50 (2012), pp. 2216-2253.

[40] V. M. VELIOV, *Second order discrete approximations to linear differential inclusions*, SIAM J. Numer. Anal. 29 (1992), pp.439-451.

[41] V. M. VELIOV). *Approximations to Differential Inclusions by Discrete Inclusions.* IIASA Working Paper. IIASA, Laxenburg, Austria: WP-89-017, 1989.